\documentclass{amsart} 
\usepackage{amssymb, amscd}
\numberwithin{equation}{section}

\def\AA{{\mathbb A}}
\def\CC{{\mathbb C}}

\def\EE{{\mathbb E}}

\def\PP{{\mathbb P}}
\def\QQ{{\mathbb Q}}

\def\ZZ{{\mathbb Z}} 

\def\k{{\mathbf k}}

\def\Acal{{\mathcal A}}

\def\Ecal{{\mathcal E}} 
\def\Fcal{{\mathcal F}}

\def\Mcal{{\mathcal M}}

\def\Ocal{{\mathcal O}}
\def\Pcal{{\mathcal P}}

\def\Xcal{{\mathcal X}}

\newcommand\Spec{\operatorname{Spec}}
\newcommand\Fil{\operatorname{Fil}}

\newcommand\Tr{\operatorname{Tr}}
\providecommand\Sp{}
\renewcommand\Sp{\operatorname{Sp}}

\let\co\colon
\newcommand\Definition[1]{\emph{#1}}
\newcommand\map[3]{\ensuremath{{#1}\co{#2}\to{#3}}}

\newcommand\proofsquare{\nobreak\hfill \hbox{%
\vrule height 5pt 
\kern-.4pt
 \vbox{%
\hrule width 5pt depth0pt height.4pt
 \kern4.6pt \hrule  }
\kern-3.75pt 
\vrule height 5pt}\kern1pt
\par}

\newtheorem{theorem}{Theorem}[section]
\newtheorem{lemma}[theorem]{Lemma}
\newtheorem{proposition}[theorem]{Proposition}
\newtheorem{corollary}[theorem]{Corollary}

\newtheorem{definition-lemma}[theorem]{Definition-Lemma}

\theoremstyle{definition}

\newtheorem{example}[theorem]{Example}

\theoremstyle{remark} 
\newtheorem{remark}[theorem]{Remark}

\begin{document}

\title[The Order of the Top Chern class of the Hodge bundle]{The Order of
the Top Chern Class of the \\ Hodge Bundle on  
the Moduli Space \\ of Abelian Varieties
\footnote{\tt topchernorder.tex\today}} 
\author{Torsten Ekedahl}
\address{Department of Mathematics\\
 Stockholm University\\
 S-106 91  Stockholm\\
Sweden}
\email{teke@math.su.se}
\author{Gerard van der Geer}

\address{Faculteit Wiskunde en Informatica, University of
Amsterdam, Plantage Muidergracht 24, 1018 TV Amsterdam, The Netherlands}

\email{geer@science.uva.nl} 


\begin{abstract}

We give upper and lower bounds for the order of the top Chern class of the Hodge
bundle on the moduli space of principally polarized abelian varieties.  

\end{abstract}

\maketitle

\begin{section}{Introduction}
\label{sec: intro}
\bigskip
\noindent
Let ${\Acal}_g/ \ZZ$ denote the moduli stack of principally polarized abelian
varieties of dimension~$g$. This is an irreducible algebraic stack of relative
dimension $g(g+1)/2$ with irreducible fibres over $\ZZ$. The stack ${\Acal}_g$
carries a locally free sheaf $\EE$ of rank $g$, the Hodge bundle, defined as
follows.  If $A/S$ is an abelian scheme over $S$ with zero section $s$ we get a
locally free sheaf $s^*\Omega^1_{A/S}$ of rank $g$ on $S$ and this is compatible
with pull backs. If $\pi: A \to S$ denotes the structure map it satisfies the
property $\Omega^1_{ A/S}= \pi^*(\EE)$ and we will consider its Chern classes
$\lambda_i(A/S):=c_i(\Omega^1_{ A/S})$ (in the Chow ring of $S$). These then are
the pullbacks of the corresponding classes in the universal case $\lambda_i :=
c_i(\EE)$. The Hodge bundle can be extended to a locally free sheaf (again
denoted by) $\EE$ on every smooth toroidal compactification $\tilde {\Acal}_g$
of ${\Acal}_g$ of the type constructed in \cite{F-C}, see Ch.\ VI,4 there. By a
slight abuse of notation we will continue to use the notation $\lambda_i$ for
its Chern classes.

The classes $\lambda_i$ are defined over $\ZZ$ and give for each fibre
${\Acal}_g\otimes k$ rise to classes, also denoted $\lambda_i$, in the Chow ring
$CH^*({\Acal}_g\otimes k )$, and in $CH^*(\tilde{\Acal}_g\otimes k)$. They
generate subrings ($\QQ$-subalgebras) of $CH_{\QQ}^*({\Acal}_g\otimes k)$ and of
$CH_{\QQ}^*(\tilde {\Acal}_g\otimes k)$ which are called the
\Definition{tautological subrings}.

It was proved in \cite{vdG1} by an application of the Grothendieck-Riemann-Roch
theorem that these classes in the Chow ring $CH_{\QQ}^*( {\Acal}_g)$ with
rational coefficients satisfy the following relation
\begin{equation}\label{fund rel}
(1+\lambda_1+\ldots +\lambda_g)(1-\lambda_1+\ldots +(-1)^g
\lambda_g)=1. %\eqno(1)
\end{equation}
Furthermore, it was proved that $\lambda_g$ vanishes in the Chow group
$CH_{\QQ}({\Acal}_g)$ with rational coefficients.  The class $\lambda_g$ does
not vanish on $\tilde{\Acal}_g$.  This raises two questions.  First, since
$\lambda_g$ is a torsion class on ${\Acal}_g$ we may ask what its order
is. Second, since $\lambda_g$ up to torsion comes from a class on the `boundary'
$\tilde{\Acal}_g - {\Acal}_g$ we may ask for a description of this class. As an
answer to the first question we give an upper bound on the order of $\lambda_g$ in
the third section and a non-vanishing result in the fourth section which implies
a lower bound. That result is obtained as a consequence of a more precise result
that determines the order of the Chern classes of the de Rham bundle up to a
(multiplicative) factor two. 
As an answer to the second question 
we shall generalize the well-known relation $12\, \lambda_1=\delta$
for $g=1$ to higher $g$ in a sequel to this paper \cite{E-vdG}.

The authors thank the referee for his many useful
remarks and thank the editors for their patience.

\end{section}
\begin{section}{$K$-theory for stacks}

We will need to extend some results on $K$-theory from schemes to stacks. Just
as in the schematic case we denote by $K^0(X)$ the Grothendieck group of vector
bundles on the algebraic stack $X$, and by $K_0(X)$ the Grothendieck group of
coherent sheaves. We shall use \cite{L-MB} for the general notions concerning
algebraic stacks, in particular the definition of coherent sheaf on an algebraic
stack, cf., \cite[15.1]{L-MB}. We shall follow standard usage in that by an
\emph{algebraic stack} we mean a general algebraic stack -- what is also called
an Artin stack -- while by Deligne-Mumford stack we mean an algebraic with an
\'etale chart. The first result is the homotopy exact sequence (cf.,
\cite[IX:Prop.~1.1]{SGA6}).
\begin{remark}
Our purpose is quite restricted, we only want results that can be used in the
next section, and hence we will make no attempt at maximum generality. It will
be clear that some of the results are in fact true in larger generality than
stated.
\end{remark}
We will in this section only consider algebraic stacks of finite type over a base field
$\k$.
\begin{proposition}\label{excision}
Let $X$ be a noetherian algebraic stack, \map{i}{X'}X a closed substack of $X$, and \map
jUX the complement of $X'$. Then the sequence
\begin{displaymath}
K_0(X') \stackrel{i_*}\longrightarrow K_0(X) \stackrel{j^*}\longrightarrow
K_0(U) \to 0
\end{displaymath}
is exact.
\begin{proof}
The proof is the same as that of \cite[IX:Prop.~1.1]{SGA6} once we know
that every coherent sheaf on $U$ extends to one on $X$ and that every coherent
subsheaf of an extended quasi-coherent sheaf extends as a subsheaf. Such
extensions are provided by \cite[Cor.~15.5]{L-MB}.
\end{proof}
\end{proposition}
The next step is to prove the analogue of \cite[IX:Prop.~1.6]{SGA6}, the
calculation of $K_0$ for a vector bundle.
\begin{proposition}\label{vector bundle}
Let $X$ be a noetherian Deligne-Mumford stack and $E \to X$ a vector bundle over
$X$. Then the pullback map $K_0(X) \to K_0(E)$ is an isomorphism.
\begin{proof}
As a general remark we note that $K_0$ is contravariant not just for flat maps but
for all morphisms between noetherian stacks of finite Tor-dimension.

We first prove the proposition for the case when $E = X\times\AA^1$ and for that
we follow the proof of \cite[IX:Prop.~1.6]{SGA6}. We start by noting that the
zero section of $X\times\AA^1$ is of finite Tor-dimension and hence we have a
map $K_0(X\times\AA^1) \to K_0(X)$ such that the composite $K_0(X) \to
K_0(X\times\AA^1) \to K_0(X)$ is the identity. This shows that $K_0(X) \to
K_0(E)$ is injective and it remains to show surjectivity.

We may now by noetherian induction assume that the corresponding map is an
isomorphism for every proper closed substack. By \ref{excision} it is hence
enough to show the result for some non-empty open subset of $X$. By
\cite[Prop.~6.1.1]{L-MB} there is a non-empty open subset of $X$ which is the
stack quotient of an action of a finite group on an affine scheme and  we
consider that open subset so that we may assume that $X$ itself is such a
quotient. If the affine scheme is $\Spec R$ and the group is $G$ then the
category of coherent $\Ocal_X$-modules is equivalent to the category of
finitely generated
modules over the twisted group ring $R[G]$ (``twisted'' because of the relation
$g\lambda=\lambda^gg$ for $\lambda \in R$ and $g \in G$) and the category of
coherent modules over $X \times \AA^1$ is then equivalent to the category of
finitely generated modules over $R[G][T]$, the polynomial ring over $R[G]$.
Note also that if $S$ is the ring of invariants of $G$ on $R$, then $S$ is
noetherian as $R$ is and $R[G]$ is a finite $S$-algebra. We then conclude by
\cite[XII:Th.~4.1]{Ba} which says exactly that for a finite $S$-algebra $A$, the
map $G_0(A) \to G_0(A[T])$ is an isomorphism; $G_0(B)$ in Bass' notation being
the Grothendieck group of finitely generated $B$-modules.

In the general case, we will say that two maps \map{f,g}YZ between algebraic
stacks are \emph{affine-homotopic} if there is a map \map{F}{Y\times\AA^1}Z
which restricted to the zero-section is $f$ and restricted to the $1$-section is
$g$. Note also that the $0$- and $1$-section in $X\times\AA^1$ are of finite
Tor-dimension which means that we may pull back elements of $K_0(Y\times\AA^1)$
along them. By the special case of $E=X\times\AA^1$ just proven we get that
$0^*$ and $1^*$ are both inverses to the said isomorphism and hence they are
equal. This implies that if we have two maps \map{f,g}YZ of finite Tor-dimension
related by an affine homotopy of finite Tor-dimension we get that $f^*=g^*$ on
$K_0(X)$.

This can now be applied to the identity map and the composite of the structure
map $E \to X$ and the zero-section $X \to E$ which are affine-homotopic
by the usual map $(v,t) \mapsto tv$.
\end{proof}
\end{proposition}
\begin{remark}
In \cite[IX:Prop.~1.6]{SGA6} noetherian induction is used to allow one to assume
that the scheme is reduced. Then a passage to the limit is made so that
instead of reducing to an open affine subset one reduces to a field in which
case one is dealing with regular schemes and can switch to $K^0$ instead where
the result is easier to prove. This is not possible in our situation as even for
a field $K$ the twisted group ring $K[G]$ is not of finite global dimension when
the characteristic divides the order of the subgroup of elements of $G$ acting
trivially on $K$. In \cite{Ba} the same strategy as in \cite[IX:Prop.~1.6]{SGA6}
is used to reduce to the case of a finite dimensional algebra over a field. There
a further reduction is made by dividing out by the radical of the algebra, a
step which corresponds to assuming that the scheme is reduced.
\end{remark}
We will now need to consider the topological filtration on $K$-theory. Hence we
define, for $X$ a regular algebraic stack, $\Fil^{\ge i} \subseteq K_0(X)$ to be
the subgroup generated by classes of sheaves with support of codimension $\ge
i$. Using (\ref{vector bundle}) we get an isomorphism $K_0(X) \to K_0(E)$ for
every vector bundle $E \to X$. We now define $\Fil_b^{\ge i}(X)$ to be the limit
of $\Fil^{\ge i}(E)$ for all vector bundles $E \to X$ and surjective vector
bundle maps between them.
\begin{proposition}\label{Chern classes and multiplicativity}
Let $X$ be a Deligne-Mumford stack which is the stack quotient of the action of
a finite group $G$ on a smooth quasi-projective $\k$-scheme $Y$.

i) The forgetful map $K^0(X) \to K_0(X)$ is an isomorphism.

ii) If $\dim X=n$ and $e \in K^0(X)$ maps to $\Fil_b^{\ge i}(X)$ then $c_j(e)\cap[X]
\in A_{n-j}(X)$ is zero for $0<j<i$.

iii) Using the isomorphism $K^0(X) \to K_0(X)$ to get a multiplication on
$K_0(X)$ the filtration $\Fil_b^\cdot$ is multiplicative, i.e.,
$\Fil_b^i\cdot\Fil_b^j\subseteq \Fil_b^{i+j}$ for all $i$ and $j$.
\begin{proof}
For i) we use the fact that the category of coherent $\Ocal_X$-modules is
equivalent to the category of coherent $\Ocal_Y$-modules with a $G$-action
compatible with the action of $G$ on $Y$, i.e., a $\Ocal_Y[G]$-module for the
twisted group ring. What needs to be proved is, as in the spatial (i.e.,
algebraic space) case, that every coherent $\Ocal_X$-module has a finite
resolution by coherent locally free $\Ocal_X$-modules. An $\Ocal_X$-module is
locally free exactly when the corresponding $\Ocal_Y$-module is. As $Y$ is
regular and of finite type over $\k$ there is a bound for the global dimension
of the local rings of $Y$. Hence it is enough to show that each coherent
$\Ocal_Y[G]$-module $\Fcal$ is the quotient of a $\Ocal_Y[G]$-module $\Ecal$
that is coherent and locally free as $\Ocal_Y$-module. By
\cite[Cor. II:2.2.7.1]{SGA6}, $\Fcal$ is the quotient of a coherent locally free
$\Ocal_Y$-module $\Ecal'$ and it is then the quotient of the $\Ocal_Y[G]$-module
$\Ecal'[G]$.

As for ii) it follows directly from i), excision (\cite[Prop. 2.3.6]{Kr}),
boundedness by dimension (\cite[Prop. 3.4.2]{Kr}), and the fact that vector
bundle maps induce isomorphisms on Chow groups (\cite[Cor. 2.4.9]{Kr}).

Continuing with iii), we may, after possibly replacing $X$ by a vector bundle
over it, reduce to showing that $\Fil^i\cdot\Fil^j\subseteq \Fil_b^{i+j}$.
According to \cite[Prop. 3.5.6]{Kr} we may find a vector bundle $E \to X$ and an
open subset $U \subseteq E$ whose complement has codimension $> i+j$ such that
$U$ is an algebraic space. For our purposes we will however need $U$ to be
quasi-projective. This is easily arranged by letting $\pi\colon E \to X$ be the
stack quotient of $G$ acting on $\Ocal_Y[G]^n$ for $n$ large ($n > i+j$ will do)
and letting $U$ be the quotient of the part where $G$ acts freely. We now have
an isomorphism $\pi^*\colon K_0(X) \to K_0(E)$ and it will be enough to show
that $\pi^*\Fil^i(X)\cdot\pi^*\Fil^j(X)\subseteq \Fil^{i+j}(E)$ and for that is
enough that $\Fil^i(E)\cdot\Fil^j(E)\subseteq \Fil^{i+j}(E)$. By excision,
(\ref{excision}), it is enough to show that $\Fil^i(U)\cdot\Fil^j(U)\subseteq
\Fil^{i+j}(U)$ but $U$ is a smooth quasi-projective variety and then it is
\cite[0 App:Cor. 2.12.1]{SGA6} when the base field is algebraically closed. The
only dependence on the assumption of algebraic closedness is for the moving
lemma. However, the moving lemma is true over any field, cf., \cite[\S3:Thm]{R}.
\end{proof}
\end{proposition}
\begin{remark}
We have used the notation $\Fil$ for the filtration directly defined by a
support condition as it is already well-established in the spatial
case. However, we do not believe it is the ``right'' definition for a general
algebraic stack. Consider for instance the case of the stack quotient, $[*/G]$,
of a finite group $G$ acting on $*=\Spec\k$. In that case $K_0([*/G])$ equals
the representation ring of $G$-representation over $\k$ and
$\Fil^{1}(*/G)=\{0\}$ which seems to be too small. On the other hand, consider
the case when $G=\ZZ/2$ (and $\k$ is a field of characteristic different from
$2$). For any $G$-representation $V$ on which $G$ acts non-trivially we have
that $\Fil_b^{\ge 1}K_0([V/G])$ equals the group $I$ of virtual bundles of rank
$0$ (which is what one would like). Using the multiplicativity we get that
$\Fil_b^{\ge n} \supseteq I^n$. On the other hand, using (\ref{Chern classes and
multiplicativity}:ii) one can show the opposite inclusion and hence $\Fil_b^{\ge
n} = I^n$. This equals the topological filtration on $K^0(BG)$ making it seem
quite reasonable.
\end{remark}
As has been stated in the introduction one of our goals is to give an explicit
integer killing $\lambda_g$ on $\Acal_g$. To make such a statement as useful as
possible one would like to be able to conclude from this that its pullback along
any map $S \to \Acal_g$ vanishes. In the case of schemes such a conclusion is
possible because the Chern classes lie in Chow cohomology groups which is a
contravariant functor. It is of course possible to formulate this contravariant
character of Chern classes without introducing Chow cohomology groups but it
would be quite awkward particularly when it comes to expressing relations
between them. We shall therefore very briefly introduce, by perfect analogy with
the spatial case, Chow cohomology groups. We shall only prove the minimal
results necessary to formulate and prove our result on $\lambda_g$. Recall
(\cite{Kr}) that an algebraic stack is said to be \emph{filtered by global
quotients} if it has a stratification by substacks such that each of them is the
quotient of an algebraic space by an affine group scheme. Note that the refined
Gysin maps (which will be used in our proof) currently (cf., \cite[5.1]{Kr}) are
defined (essentially) only in the case when the involved stacks are filtered by
global quotients. We now follow \cite{Fu} in defining for a map $X \to Y$ of
algebraic stacks filtered by global quotients the \emph{bivariant Chow groups}
$A^p(X \to Y)$ consisting of collections of operations as in \cite[Def.\
17.1]{Fu} fulfilling the conditions $C1-C3$ with the difference that in $C3$ the
map $i\colon Z'' \to Z'$ is assumed to be representable, locally separated whose
normal cone stack is a vector bundle stack of constant rank. (These are on the
one hand the conditions under which Kresch defines the refined Gysin map, on the
other hand we shall want to use the Gysin map for the diagonal of an algebraic
stack which means that it will not be enough to require $C3$ for
l.c.i. embeddings.) As in [loc.\ cit.] we define the Chow cohomology groups to
be the bivariant groups for the identity maps.
\begin{proposition}\label{Chow cohomology}
i) Let $X$ be an algebraic stack filtered by global quotients and $E$ a vector bundle
over $X$. If for a map $g\colon X' \to X$ we set $c_p(E)(\alpha):=c_p(g^*E)\cap
\alpha$ then we obtain an element of $A^p(X)$.

ii) Let $X$ be a smooth algebraic stack of pure dimension $n$ that is filtered
by global quotients. The map $A^p(X) \to A_{n-p}(X)$ given by cupping an
operation with the fundamental class is an isomorphism.
\begin{proof}
We will follow very closely the relevant parts of \cite[Chap.\ 17]{Fu}. For the
first part we need the compatibility of Chern classes with proper push forwards,
pullbacks and Gysin maps respectively. The first two properties are proven in
\cite[\S2.5]{Kr} and the proof of the third follows quite directly:\footnote{We
are grateful to Andrew Kresch for providing us with the following proof.} Chern
classes are polynomials in Segre classes so it is enough to prove that Segre
classes commute with Gysin maps. From the definition of Segre classes
(\cite[Def.\ 2.5.4]{Kr}) we see that they in turn are expressed in terms of flat
pullbacks, (iterated) top Chern class operations, and proper push
downs. Commutation of Gysin maps with flat pullbacks and proper push downs are
clear (cf., \cite[\S5.1]{Kr}) so it remains to show that they commute with top
Chern classes. Looking at the construction of the Gysin map (cf.,
\cite[\S5.1,\S3.1]{Kr}) one sees that it is expressed in terms of flat
pullbacks, proper push forwards, and intersection with a principal effective
Cartier divisor (the normal cone $C_{F'}G'$ in the deformation space
$M^o_{F'}G'$, cf., \cite[\S5.1]{Fu}). As the other operations are already known
to commute with the top Chern class one is reduced to proving commutation with
intersection with a principal effective Cartier divisor (the operation
introduced in \cite[\S2.2]{Kr}). Looking at the definition of the top Chern
class what is needed is the projection formula (cf., \cite[Prop.~2.3 c]{Fu}) for
intersection with a divisor. This is proved as in the proof of
\cite[Prop.~2.3]{Fu}.

As for ii), the proof is identical with \cite[Thm.\ 17.4.2]{Fu} (specialized to
the case of the base being a point). (Note that in that proof we can not assume
that the diagonal map is an immersion which is why we have to allow for more
general maps in $C3$.)
\end{proof}
\end{proposition}
\end{section}
\begin{section}{A bound on the order of the class $\lambda_g$.}

We will make some computations in the Chow group of $\Acal_g$. To make this
reference applicable (and for other reasons) all our algebraic stacks will be of
finite type over a field.  We begin with a lemma which is no doubt well known.
\begin{lemma}\label{chern}
If $E$ is a vector bundle of rank $g$, the total Chern class of the
graded vector bundle $\Lambda^*E$ is zero in degrees 1 to $g-1$ and
$-(g-1)!c_g(E)$ in degree $g$.
\end{lemma}
%% \noindent
%% {\sl Proof.} 
\begin{proof}
As usual we note that the components of the total Chern class are universal
polynomials with integer coefficients in the Chern classes $c_i:=c_i(E)$. Such a
relation is true if and only if it is true for the tautological bundle on all
Grassmannians of $g$-dimensional subspaces and we may hence let $E$ be such a
bundle. Then as we may think of the coefficients in the universal polynomials as
rational numbers we can note that the Chern classes of degree 1 to $g-1$ of
$\Lambda^*E$ vanish if and only if the Newton polynomials of the same degrees
do, that is, if and only if ${\text ch}(\Lambda^*E)$ vanishes in the same
degrees. Thus the first part follows from the Borel-Serre formula
(cf.~\cite{B-S})
$$
{\text ch}(\Lambda^*E)=(-1)^gc_g(E){\text{Td}}(E)^{-1}.
$$
Furthermore, if the Chern classes of degree 1 to $g-1$ vanish for a
bundle $F$, then it is clear that $(-1)^{g-1}g\, c_g(F)=s_g(F)$, as can be
seen from Newton's formula. (Here $s_i(F)$ are the Newton polynomials in the
(roots of the) Chern classes of $F$.) Hence in degree $g$ we have  
${\text ch}(\Lambda^*E)=
(-1)^{g-1}gc_g(\Lambda^*E)/g!$, and using the Borel-Serre formula again gives
the desired formula.
%%%$
\end{proof}
\begin{lemma}\label{orderp}
Let $p$ be a prime. Suppose that $\pi\colon A \to S$ is a family of abelian varieties
of relative dimension $g$, where $S$ is an algebraic stack that is the quotient
of a smooth quasi-projective $\k$-scheme by a finite group, and assume that $L$
is a line bundle on $A$ of order $p$, on all fibres of $\pi$. Let $E$ be the
Hodge bundle of $\pi$, i.e., the pullback of $\Omega^1_{A/S}$ along the zero
section and let $e$ be the class in $K^0(S)$ of the graded bundle
$\Lambda^*E$. If $p>2g$ then $pe \in \Fil_b^{\ge g+1}$ and in particular
$p(g-1)!\lambda_g=0$. If $S$ is actually a smooth quasi-projective variety we
have that the stronger conclusion $pe \in \Fil^{\ge g+1}$ is true provided only
that $p > \min(2g,\dim S+g)$.
\end{lemma}
%% $
%% {\sl Proof.}  
\begin{proof}
That $pe \in \Fil_b^{\ge g+1}$ implies that $p(g-1)!\lambda_g=0$ follows from
(\ref{chern}), (\ref{Chern classes and multiplicativity}:ii), (\ref{Chow
cohomology}). 

As in the proof of (\ref{Chern classes and multiplicativity}) we can find a
vector bundle $E \to S$ with an open set $U$ whose complement has codimension
$>g$ that is a smooth quasi-projective variety. By (\ref{vector bundle}) and the
definition of $\Fil_b$ we may pull back $\pi$ to $U$ and hence assume that $S$
is smooth and quasi-projective.

By twisting $L$ by a line bundle on $S$ so that it is trivial along the zero
section we may assume that it is of order $p$ on $A$.  Denoting the class of $L$
in $K_0$ by $[L]$ we then either have that $p([L]-1)$ has support of codimension
$>2g$ if $p>2g$ or is zero if $p>\dim S+g$. Indeed, from the relation
$$
0=[L]^p-1=p([L]-1)(1+\frac{p-1}{2}([L]-1)+...)+([L]-1)^p
$$ 
and the fact that $[L]-1$ is nilpotent (which in turn follows for instance from
the multiplicativity of the topological filtration) we get that $p([L]-1)$ is
$([L]-1)^p$ times a unit and $([L]-1)^p$ element is supported in codimension
$\ge p$ by the multiplicativity of the topological filtration. Now in the first
case the image under $\pi$ of the support has codimension $>g$ on $S$ and so we
may safely remove it and may assume that $p[L]=p$ in $K_0(A)$.

Consider now the Poincar\'e bundle ${\Pcal}$ on $A\times_S\check A$, where
$\check A$ denotes the dual abelian variety.  By base change $R\pi_*{\Ocal}_A$
is the (derived) pullback along the zero-section of $\check A$ of the sheaf
$R\pi_*\Pcal$. We have that $p[{\Pcal}]=p[L\bigotimes{\Pcal}]$ and so
$p[R\pi_*{\Pcal}]=p[R\pi_*(L\bigotimes {\Pcal})]$. Now, a fibrewise calculation
shows that $R\pi_*(L\bigotimes \Pcal)$ has support along the inverse of the
section of $\check A$ corresponding to $L$. As that section is everywhere
disjoint from the zero section the pullback of $R\pi_*(L\bigotimes \Pcal)$ along
the zero section is 0 and thus $p[R\pi_*{\Ocal}_A]=0$. Now, the well-known
calculation of  $R^i\pi_*{\Ocal}_A$ shows that it is isomorphic to $\Lambda^iE$
so the lemma follows.
\end{proof}
\begin{definition-lemma}\label{lcd}
For an integer $g$ we let $n_g$ be the largest common divisor of all $p^{2g}-1$
where $p$ runs through all primes larger than a sufficiently large fixed number
$N$ (which may be taken to be $2g+1)$). For an odd prime $\ell$, the exponent
$k$ of the exact power $\ell^k$ of $\ell$ that divides $n_g$ is the largest $k$
such that $\ell^{k-1}(\ell-1)$ divides $2g$ and $0$ if $\ell-1$ does not divide
$2g$. The exponent $k$ of the exact power $2^k$ that divides $n_g$ is the
largest $k$ such that $2^{k-2}$ divides $2g$.
\end{definition-lemma}
\begin{proof}
Note that $\ell^k$ divides $n_g$ if and only if the exponent of 
$(\ZZ/\ell^k\ZZ)^*$ divides $2g$.
The statement now follows directly from the structure of $(\ZZ/{\ell}^k)^*$ and Dirichlet's prime
number theorem.
\end{proof}
\smallskip
\begin{example}
We have $n_1=24$, $n_2=240$, $n_3=504$ and $n_4=480$.
\end{example}
\begin{theorem}\label{lambdascheme}
Suppose that $\pi\colon A \to \Acal_g$ is the universal family of principally
polarized abelian varieties of relative dimension $g$. Then
$(g-1)!n_g\lambda_g=0$ on $\Acal_g$.
\end{theorem}
%% \noindent
%% {\sl Proof.} 
\begin{proof}
For any sufficiently large prime $p$ we can apply lemma \ref{orderp} on the
cover obtained by adding a line bundle everywhere of order $p$. Projecting down
to $\Acal_g$ again and using that that cover has degree $p^{2g}-1$ (being equal
to the number of line bundles of order $p$) gives
$(g-1)!p(p^{2g}-1)\lambda_g=0$. We then finish by using Definition \ref{lcd}
(and noting that the factor $p$ causes no trouble as by using several primes we
see that no prime $>2g$ can divide the smallest annihilating integer).
\end{proof}
\begin{remark}
For various mostly technical reasons we have been working over a field,
where of course the prime fields are the optimal choices. If one would like to
prove the statement over the integers then assuming that the technical details can be
overcome we will still have to deal with the fact that one needs to invert a
finite number of primes to use the trick of producing a line bundle of prime
order by going to a covering. Though there is a large freedom in
choosing the primes to which we could apply this trick, still we do not know how
to solve this problem.
\end{remark}
\end{section}
\begin{section}{The order of the Chern Classes of the de Rham Bundle}

We will now consider the Chern classes of the bundle of first relative de Rham
cohomology of the universal abelian variety over $\Acal_g$. Over the complex
numbers these Chern classes are the Chern classes of a flat bundle. By the
Chern-Weil expression of Chern classes in terms of curvature for a smooth
manifold, the Chern classes are torsion classes in integral
cohomology. Grothendieck has given an arithmetic proof of this fact. His proof
gives an explicit bound for the orders of the Chern classes which we will
exploit. This bound depends on the field of definition of the representation of
the fundamental group that is associated to the flat bundle. Luckily, in our
case this representation is the natural representation of $\Sp_{2g}(\ZZ)$
and hence the field of definition is $\QQ$ which makes the bounds the best
possible.

We will then give lower bounds for the order of these Chern classes. From the
complex point of view, where the cohomology of $\Acal_g$ can be thought of as
$H^*(\Sp_{2g}(\ZZ),\ZZ)$, the idea is to find finite (cyclic) subgroups of
$\Sp_{2g}(\ZZ)$ and compute the order of the restriction of the Chern classes to
the cohomology of such a subgroup. Now, any finite subgroup $G$ of
$\Sp_{2g}(\ZZ)$ is actually the group of automorphisms of a principally
polarized variety $A$ (as such a subgroup has a fixed point on the Siegel upper
half space). Such an $A$ gives not only a point of $\Acal_g$ but also a map from
the stack quotient $[*/G]$ to $\Acal_g$. In particular it gives a map from the
cohomology of $\Acal_g$ to that of $G$. This is an algebraic construction of the
map given by the map from $G$ to $\Sp_{2g}(\ZZ)$. Even though it would be
possible to give a purely arithmetic construction of the detecting subgroups
that we will consider it will be (at least) as easy to construct the principal
abelian variety with a group action and we shall do exactly that. We shall then
also continue to use the stack language. An extra benefit of this way of
presenting the argument is that we directly get the lower bound also in positive
characteristic. In that case, if one knew that the specialisation map in the
cohomology of $\Acal_g$ induced an isomorphism then the result would follow from
the complex version. Though such a specialisation result seems rather
straightforward using the proper and smooth base change theorems, the toroidal
compactifications of Chai and Faltings and an induction on $g$ we know of no
reference. Our argument will avoid reference to such a specialisation result.

Even though we shall apply our results to obtain information on the order of
$\lambda_g$, the results we obtain should be of independent interest. 

We start with some preliminary comments on the cohomology of stacks and in
particular how to extend results from the cohomology of schemes to that of
algebraic stacks. If $\Acal \to S$ is an algebraic stack over a scheme $S$ we
may find a chart $U \to \Acal$, i.e., a smooth surjective map from a scheme
$U$. We may then consider the simplicial $S$-scheme ${U\times_{\Acal}}_\bullet$
which in degree $n$ is $U\times_{\Acal}U\times_{\Acal}U\dots U\times_{\Acal}U$
($n+1$ times) with the obvious face and degeneracy maps. We then have an
augmentation map ${U\times_{\Acal}}_\bullet \to \Acal$ which by smooth descent
induces an isomorphism in cohomology. For any simplicial scheme $X_\bullet$ there
is (cf., \cite[Exp $V^{bis}$]{SGA4}) a spectral sequence which at the $E_1$-term
is the cohomology of $X_n$ (with coefficients in the sheaf whose cohomology
should be computed) converging to the cohomology of $X_\bullet$. We may therefore
reduce a number of questions on the cohomology of stacks to that of schemes:
\begin{itemize}
\item If $S = \Spec\CC$ we may use the comparison theorem between \'etale and
classical cohomology with finite (or $\ell$-adic) coefficients to get a
comparison theorem for the \'etale and classical cohomology of algebraic stacks.

\item If $S$ is a discrete valuation ring then for any smooth map $X \to S$ the
smooth base change theorem gives us a specialisation map (cf., \cite[Exp.\
XIII]{SGA7} or \cite[V, (1.6.1)]{Arcata} where it is called the ``cospecialisation
map'') $H^*(X_{\overline\eta},\ZZ/n) \to H^*(X_{\overline s},\ZZ/n)$, where $n$
is invertible on $S$, $\overline s$ is a geometric point over the special fibre
of $S$, and $\overline\eta$ a geometric point over the generic. This extends
immediately to smooth simplicial $S$-schemes and hence to smooth $S$-stacks.

\item If $S = \Spec\CC$ we get that the de Rham Chern classes of a flat vector
bundle over a smooth stack are trivial. In fact, the Chern-Weil construction has
been done for simplicial manifolds (cf., \cite{Du}).
\end{itemize}
When we will be talking about $\ell$-adic cohomology for an algebraic stack
\map{\pi}{\Acal}{\Spec\k} over a field we shall always mean the higher direct
images of $\pi$, or equivalently the cohomology of the pullback of $\Acal$ to a
separable closure of $\k$.

Note now that the \Definition{de Rham bundle}, $H_{dR}^1:=
R^1\pi_*\Omega^\cdot_{{\Xcal}_g/{\Acal}_g}$, where \map\pi{{\Xcal}_g}{{\Acal}_g}
is the universal family, is indeed provided with an integrable connection and
hence its Chern classes in integral and hence, by the comparison and
specialisation theorems, $\ell$-adic cohomology\footnote{$\ell$ being as usual
a prime different from the characteristic.}, which we will denote $r_i \in
H^{2i}(\Acal_g,\ZZ_\ell(i))$, are torsion classes. In this section we shall
determine their exact order up to a factor of $2$.

We begin by using a result of Grothendieck to get an upper bound for the order
of $r_i$.
%;;%%%Upper bound for order Proposition upper bound
\begin{proposition}\label{upper bound}
i) We have $r_i=0$ for odd $i$.

ii) We have $n_ir_{2i}=0$.
\end{proposition}
\begin{proof}
The first part follows immediately because $H_{dR}^1$ is a symplectic vector
bundle.

As for the second part we may assume that the characteristic is $0$ as the case
of positive characteristic follows from the characteristic $0$ case by
specialisation as the Chern classes are compatible with specialisation maps. In
that case we may further reduce to the case of the base field being the complex
numbers. We may also prove the annihilation of $n_ir_{2i}$ in $\ell$-adic
cohomology for a specific (but arbitrary) prime $\ell$.

Now, the existence of the Gauss-Manin connection on $H^1_{dR}$ means that
$H^1_{dR}$ has a discrete structure group. More precisely, the (classical)
fundamental group of the algebraic stack ${\Acal}_g$ is $\Sp_{2g}(\ZZ)$ and
$H^1_{dR}$ is the vector bundle associated to the representation of it given by
the natural inclusion of $\Sp_{2g}(\ZZ)$ in $\Sp_{2g}(\CC)$. This complex
representation is obviously defined over the rational numbers so we may apply
\cite[4.8]{Gr} with field of definition $\QQ$. We thus conclude that $r_{2i} \in
H^{2i}({\Acal}_g,\ZZ_\ell(i))$ (for the analytic, and hence by the comparison
theorem, the algebraic stack) is killed by $\ell^{\alpha(i)}$, where $\alpha(i)$
is defined as
\begin{displaymath}
\inf_{\lambda \in H}v_\ell(\lambda^i-1)
\end{displaymath}
and $H \subseteq \ZZ_{\ell}^\ast$ is the image of the Galois group of the field
of definition of the cyclotomic character. However, as the base field is $\QQ$
this image is all of $\ZZ_{\ell}^\ast$ and the result follows from the
definition of $n_i$.
\end{proof}
We now aim to show that this upper bound is the precise order of the $r_i$,
using as was explained actions of finite cyclic groups on principally polarized
abelian varieties.
%;;%%%Lower bound for order Proposition lower bound
\begin{proposition}\label{lower bound}
Assume that $i \le g$. The order of $r_{2i}$ is divisible by $n_i/2$ over
$\CC$. In general, for each prime $\ell$ different from the characteristic of the
base field, $r_{2i}$ in $\ell$-adic cohomology has order divisible by the
$\ell$-part of $n_i/2$.
\end{proposition}
\begin{proof}
The $\ell$-adic part implies the integral cohomology part so we may pick a prime
$\ell$ different from the characteristic of the base field and look at $r_{2i}$
in $\ell$-adic cohomology. What we want to show is that if $\ell$ is odd and
$\ell-1|2i$ and $k$ is the largest integer such that $\ell^{k-1}(\ell-1)|2i$, then
$r_{2i}$ has order at least $\ell^k$ and similarly for $\ell=2$. We will do this
by defining a map $B\ZZ/\ell^k \to \Acal_g$, such that inverse image of the
$r_{2i}$ to $H^{4i}(\ZZ/\ell^k,\ZZ_\ell)$ has order $\ell^k$. Now,
$H^{4i}(\ZZ/\ell^k,\ZZ_\ell)$ is isomorphic to $\ZZ/\ell^k$ and hence an element
in it has order $\ell^k$ if and only if its reduction modulo $\ell$ is
non-zero. As $H^{4i}(\ZZ/\ell^k,\ZZ_\ell)/\ell$ injects into
$H^{4i}(\ZZ/\ell^k,\ZZ/\ell)$ it is enough to show that the pullback of $r_{2i}$ is
non-zero in $H^{4i}(\ZZ/\ell^k,\ZZ/\ell)$. Note that a map $B\ZZ/\ell^k \to
\Acal_g$ consists of a principally polarised $g$-dimensional abelian variety $A$
over the base field together with an action on it (preserving the polarisation)
by $\ZZ/\ell^k$. The pullback of $r_{2i}$ is then obtained as follows: The group
$\ZZ/\ell^k$ acts on $H^1_{DR}(A)$. A representation of $\ZZ/\ell^k$ has Chern
classes\footnote{In the classical case only for a complex representation but
Grothendieck (cf., \cite{Gr}) extended it to a representation over any field. In
any case our representations will be ordinary and can thus be lifted to
characteristic zero.} in $H^*(\ZZ/\ell^k,\ZZ_{\ell})$ and the pullback of
$r_{2i}$ is the $2i$'th Chern class of $H^1_{DR}(A)$.

Assume now to begin with that $\ell$ is odd. Consider any Galois cover $C
\to \PP^1$ with Galois group $\ZZ/\ell^k$ which is ramified at $0$, $1$ and
$\infty$ with ramification group of order $\ell^k$, $\ell^k$, and $\ell$
respectively (the existence of such a cover follows directly from Kummer
theory). By the Hurwitz formula the genus of such a covering fulfills the
relation $2g-2=-2\ell^k+2(\ell^k-1)+\ell^{k-1}(\ell-1)$, i.e.,
$2g=\ell^{k-1}(\ell-1)$.

We claim that the action of $\ZZ/\ell^k$ on $H^1_{DR}(C)$ is isomorphic to the
sum of all primitive\footnote{I.e., of order exactly $\ell^k$.} characters of
$\ZZ/\ell^k$. Admitting that for the moment we can go on with computing its
total Chern class.  Choose a primitive $\ell^k$'th root of unity $\zeta$ and use
it in particular to identify $\ZZ/\ell$ with $\mu_\ell$. Let $x \in
H^2(\ZZ/\ell^k,\ZZ/\ell)$ be the generator given as the first Chern class of the
character that takes $1$ to $\zeta$. Then the first Chern class of the character
of $\ZZ/\ell^k$ that takes $1$ to $\zeta^i$ equals $ix$ and hence the total
Chern class of that character is $1+ix$. By the multiplicativity of the total
Chern class we get that the total Chern class of $H^1_{DR}(C)$ equals
$\prod_{(i,\ell)=1}(1+ix)$ which in turn is equal to
$(1+x^{\ell-1})^{\ell^{k-1}}=1+x^{\ell^{k-1}(\ell-1)}$ and hence the
$\ell^{k-1}(\ell-1)$'th Chern class is non-zero. If $2i$ instead is a proper
multiple $2i=m\ell^{k-1}(\ell-1)$ we look at the $m$'th power of the
Jacobian of $C$ (with the diagonal action of $\ZZ/\ell^k$).  This gives the
non-triviality when $g=i$ and when $g > i$ we simply add a principally polarised
factor on which $\ZZ/\ell^k$ acts trivially.

It remains to show that the action of $\ZZ/\ell^k$ on $H^1_{DR}(C)$ is as is
claimed. We start with a remark on Chern classes of representations of a finite
group. They only depend on the corresponding element in the representation ring
and those elements in turn are determined by their character (i.e., the traces
of the actions of the group elements). We shall therefore speak of the Chern
classes of a character.

The Lefschetz fixed point formula gives a formula for the character in terms of
the number of fixed points\footnote{Note that as the orders of these elements
are prime to the characteristic, each fixed point is counted with multiplicity
$1$.} of the elements of $\ZZ/\ell^k$. A proof using this formula is certainly
possible but the following argument requires less computations. As the Lefschetz
fixed point formula is also valid for $\ell$-adic cohomology the looked for
character is the same as that for the action on $H^1(C,\QQ_\ell)$ or more
precisely, the character of the representation on $H^1_{DR}(C)$ is the reduction
modulo $p$ of the character on $H^1(C,\QQ_\ell)$. As the action of $\ZZ/\ell^k$
on $C$ is faithful, its action on $H^1(C,\QQ_\ell)$ is faithful as well so that
at least one primitive character has to appear in $H^1(C,\QQ_\ell)$. On the
other hand the $\ell^k$'th cyclotomic polynomial is irreducible over $\QQ_\ell$
so if one primitive character appears they must all appear. Now, the number of
such characters is $\ell^{k-1}(\ell-1)$ which equals the dimension of
$H^1(C,\QQ_\ell)$ so that $H^1(C,\QQ_\ell)$ consists of each primitive character
exactly once.

When $\ell=2$ we may assume that $k > 2$ as the lower bound to be proven for
$k\leq 2$ is implied by the one for $k=3$. We then make essentially the same
construction, a Galois cover with group $\ZZ/2^{k}$ of $\PP^1$ ramified at three
points with ramification groups of order $2^k$, $2^k$, and $2$ respectively of
genus $g=2^{k-3}$. The rest of the argument is identical to the odd case.
\end{proof}
\begin{remark}
i) It follows from (\ref{fund rel}) that the $r_i$ are torsion 
already in the Chow groups. Our result gives a lower bound for this order but we 
don't know if this bound is sharp.

ii) From the complex point of view our geometric construction can be seen simply
as constructing an element of order $\ell^k$ in
$\Sp_{\ell^{k-1}(\ell-1)}(\ZZ)$. This can be done directly; in the odd
case one may consider the ring of $\ell^k$'th roots of unity $R=\ZZ[\zeta]$ with
the obvious action of $\ZZ/\ell^k$ and the symplectic form
$\langle\alpha,\overline{\beta}\rangle :=
\Tr(\alpha\beta(\zeta-\zeta^{-1})^{-\ell^{k}+\ell^{k-1}+1})$. This is obviously
a symplectic invariant form and that it is indeed an integer-valued perfect
pairing follows from the fact that the different of $R$ is the ideal generated
by $(\zeta-\zeta^{-1})^{\ell^{k}-\ell^{k-1}-1}$.

iii) When $g=\ell^k(\ell-1)/2$ we actually get a lower bound for the top Chern class
of the de Rham cohomology of the universal curve over $\Mcal_g$. However, there
is no direct analogue of the trick of adding a factor with trivial action so
this does not give a lower bound for all $g \ge \ell^k(\ell-1)/2$.
\end{remark}
\begin{theorem}
We have that $r_{2i+1}=0$ for all $i$ and that the order of $r_{2i}$ in integral
($\ell$-adic) cohomology equals (resp.\ the $\ell$-part of) $n_i/2$ or $n_i$ for
$i \le g$.
\end{theorem}
\begin{proof}
This follows immediately from Props.\ \ref{lower bound} and \ref{upper bound}.
\end{proof}
\begin{corollary} The order of $\lambda_g$ is divisible by $n_g/2$.
\end{corollary}
\begin{proof}
The top Chern class of $H_{dR}^1$ is $\lambda_g^2$.
\end{proof}
\begin{remark}
Our upper and lower bounds for $r_{2i}$ are off by a (multiplicative) factor of
$2$. Furthermore, when $g=1$ the lower bound is the correct order. It is tempting
to believe that it is the lower bound that is the correct one for all $g$ and
furthermore that that should be related to the fact that we have a symplectic rather
than a general linear representation.
\end{remark}
\end{section}
\newcommand\eprint[1]{Eprint:~\texttt{#1}}
%\eject

%
%

\end{document}